\documentclass[12pt]{article}
\usepackage{amsxtra,amssymb,amsthm,amsmath,latexsym}

\theoremstyle{plain}

\def\R{{\mathbb R}}

\def\oH{{\overset{\circ}{H}}}
\def\oH1{{\overset{\circ}{H}\kern-.02in{}^1}}

\def\bee{\begin{equation*}}
\def\eee{\end{equation*}}
\def\be{\begin{equation}}
\def\ee{\end{equation}}

\begin{document}

\title{Estimating the size of the scatterer
}

\author{Alexander G. Ramm\\
 Department  of Mathematics, Kansas State University, \\
 Manhattan, KS 66506, USA\\
ramm@ksu.edu\\
http://www.math.ksu.edu/\,$\sim $\,ramm
}

\date{}
\maketitle\thispagestyle{empty}

\begin{abstract}
\footnote{MSC: 78A45}
\footnote{Key words: scattering amplitude; size of the scatterer
 }

Formula for the size of the scatterer is derived explicitly in terms of the scattering amplitude corresponding to this scatterer. By the scatterer either a bounded obstacle $D$
or the support of the compactly supported potential is meant.

\end{abstract}

\section{Introduction}\label{S:1}

 Let $D$ be a bounded $C^2-$smooth connected domain 
 in $\R^3$, $S$ be its boundary, $N$ is the outer unit normal to $S$, $u_N$ is the normal derivative of $u$ on $S$, $S^2$
 be the unit sphere in $\R^3$, $D':=\R^3\setminus D$, 
 $k=const>0$, $\alpha, \beta\in S^2$, $v$ be the scattered
 field, $g=g(x,y):=\frac{e^{ik|x-y|}}{4\pi |x-y|}$, $\{e_j\}_{j=1}^3$ be an orthonormal basis in $\R^3$.
 We drop the $k-$dependence since $k$ is fixed. 

 Consider the scattering problem 
\be\label{e1}
\nabla^2 u+k^2u=0 \quad in \quad D', \quad u|_S=0,\quad u=u_0+v,
\ee  
where $u_0:=e^{ik\alpha\cdot x}$ is the incident field, and $v$, the scattered field, satisfies the radiation condition
\be\label{e2}
v_r-ikv=O(r^{-2}), \quad r=|x|\to \infty, \quad x/r=\beta,
\ee
uniformly with respect to $\beta$. 
Problem \eqref{e1}--\eqref{e2} is the obstacle scattering
 problem, $D$ is the obstacle, the scatterer, $u=u(x):=u(x,\alpha)$ is the scattering solution.
 
It is known (see, for example, \cite{R190}, \cite{R670}) that problem  \eqref{e1}--\eqref{e2} has a solution and this solution $u$ is unique. 

By the Green's formula one gets
\be\label{e3}
u=u_0(x,\alpha)-\int_S g(x,s)h(s)ds, \quad h:=u_N(s,\alpha).
\ee
Let $r\to \infty, x/r=\beta$ in \eqref{e3}. Then
\be\label{e4}
u=u_0+A(\beta, \alpha)\frac{e^{ikr}}{r}+O(r^{-2}),\quad r\to \infty,\quad \beta=x/r,
\ee
where 
\be\label{e5}
A(\beta, \alpha)=-\frac {1}{4\pi}\int_S e^{-ik\beta \cdot s}h(s, \alpha) ds.
\ee
The properties of $A(\beta, \alpha)$ are studied, for example, in \cite{R190}, \cite{R670}.

 It was proved in 
\cite{R190}, p.62,  that $A(\beta, \alpha)$ is an analytic function of $\beta$ and $\alpha$ in the algebraic variety  
$M\subset \mathbb{C}^3$ defined by the equation $z\cdot z=1$,
$z\in \mathbb{C}^3$, where $z\cdot z:=\sum_{j=1}^3 z_j^2$,
$z_j\in  \mathbb{C}^1$. Indeed, from \eqref{e5} it follows that $A(\beta, \alpha)$ is an analytic function of $\beta$
on the variety $M$.
This function is originally defined on $S^2$. For small $p>0$ it is defined for $|\beta_1|\le p, |\beta_2|\le p$ and
$\beta_3=(1-\beta_1^2-\beta_2^2)^{1/2}$ and is an analytic
function of the three complex variables $\beta_1, \beta_2, \beta_3$ in the region $|\beta_1|\le p, |\beta_2|\le p,
\beta_3=(1-\beta_1^2-\beta_2^2)^{1/2}$. The scattering amplitude $A(\beta, \alpha)$ admits a unique analytic continuation to the algebraic variety $M$.
Therefore, one can take $\beta$ in formula \eqref{e5}, for example,  as
\be\label{e6}
\beta=ae_1+ibe_2, \quad a,b\in \R, \quad a^2-b^2=1.
\ee
This $\beta$ clearly belongs to $M$.

Let us call a plane $P$ supporting $D$ at a point $s$ if 
$s\in D\cap P$ and $D$ is contained in one of the two half-spaces bounded by $P$. A plane is supporting $D$ if it is
supporting $D$ at some point.  The distance $d$ between two parallel planes supporting $D$ we call the size of $D$ in the direction of the normal to these planes. 
 
  Our first result is the following theorem.

{\bf Theorem 1.} {\em  The size of $D$ in the direction $e_2$
is 
\be\label{e7}
d\le\lim_{b\to \infty}\frac{\ln |A(ae_1+ibe_2, \alpha)|}{bk},
\ee
where $a^2-b^2=1$.
}

Consider now the potential scattering:
\be\label{e8}
[\nabla^2+k^2-q(x)]\psi=0\quad in \quad \R^3,
\ee
where $\psi$ satisfies the radiation condition \eqref{e2},
the potential $q$ is a real-valued compactly supported function with support $D$, $q\in L^2(D)$. It is known that 
$\psi$ satisfies relation \eqref{e4} with
\be\label{e9}
A(\beta, \alpha)=-\frac 1 {4\pi}\int_D e^{-ik\beta \cdot y}H(y)dy, \quad H(y)=H(y, \alpha):=q(y)\psi(y, \alpha).
\ee
The function $H$ does not depend on $\beta$. One has
\be\label{e9'}
\max_{y\in D, \alpha\in S^2}|H(y, \alpha)|\le c,
\ee
where $c>0$ is a constant.

Since $S$
is $C^2-$smooth, any section of $S$ by a plane is a curve with a bounded length.

Our second result is:

{\bf Theorem 2.} {\em Formula \eqref{e7} holds for the estimate of the support $D$ of $q$ in the direction $e_2$.}

In the next section proofs are given.

\section{Proofs}\label{S:2}

{\bf Proof of Theorem 1.} Using formula \eqref{e5} with
$\beta$ from \eqref{e6}, one gets
\be\label{e10}
|A(ae_1+ibe_2, \alpha)|\le c\int_S e^{bks_2}ds, \quad 
c=\frac 1{4\pi}\max_{s\in S}|u_N(s,\alpha)|.
\ee
The conclusion of Theorem 1 follows from \eqref{e10}
and Lemma 1.

{\bf Lemma 1.} {\em One has
\be\label{e11}
\lim_{b\to \infty}\frac {\ln \int_Se^{bks_2}ds}{bk}=d,
\ee
where $d$ is the size of $D$ in the direction $e_2$.} 
 
 Theorem 1 is proved. \hfill$\Box$
 
 {\bf Proof of Lemma 1.} One has
\be\label{e12}
\int_S e^{bks_2}ds=\int_{S; s_2\ge d-\epsilon}e^{bks_2}ds+
O(e^{bk(d-\epsilon)}), \quad b\to \infty.
\ee   
We assume that $S$ is $C^2-$smooth. Therefore the section of $S$ by any plane $y=const$ is a curve that has has a finite length $\ell(y)$, $\ell(y)=cy^{\gamma}$, where $\gamma>0$ depends on the smoothness of $S$ at the supporting point.
 In the coordinate system that we use one of the planes, supporting $D$, has the equation $y=0$ and the other
plane, supporting $D$, is $y=d$. If one changes variable
$y$ setting $d-y=\eta$, then
\be\label{e13}
J:=\int_{S; s_2\ge d-\epsilon}e^{bks_2}ds= \int_{d-\epsilon}^{d}e^{bky}\ell(y)dy=e^{bkd}
\int_0^{\epsilon}d\eta e^{-bk\eta}c(d-\eta)^{\gamma}.
\ee
Thus
 \be\label{e14}
 \lim_{b\to \infty}\frac {\ln J}{bk}= d.
 \ee
 Lemma 1 is proved. \hfill$\Box$
 
 {\bf Proof of Theorem 2} is essentially the same as the proof of Theorem 1.

\section{Remarks}\label{S:3}

1. We have chosen $\beta$ by formula \eqref{e6}, but we can
choose it so that Im$\beta$ is directed along a given unit vector from $S^2$.

2. The symmetry of the scatterer can influence the functional form of the scattering amplitude. 

In \cite{R691} it is
proved that if and only if the scatterer is spherically symmetric its scattering amplitude is a function of $\beta \cdot \alpha$, $A=A(\beta \cdot \alpha)$. Therefore,
when the scatterer is spherically symmetric $\alpha$
should not be chosen orthogonal to the vector that is parallel to 
Im$\beta$.

3. In quantum mechanics the scattering amplitude
corresponding to a spherically symmetric compactly supported potential is calculated as a series \cite{L}. This series is 
less convenient for passing to the limit $b\to \infty$ than 
formula \eqref{e9}.

\newpage

\end{document}